\documentclass[12pt]{article}
\def\bu{\bullet}
\def\marker{\>\hbox{${\vcenter{\vbox{
    \hrule height 0.4pt\hbox{\vrule width 0.4pt height 6pt
    \kern6pt\vrule width 0.4pt}\hrule height 0.4pt}}}$}\>}
\def\gpic#1{#1
     \medskip\par\noindent{\centerline{\box\graph}} \medskip}

\usepackage{amsmath, amssymb, amsthm, fullpage}
\usepackage{verbatim}

\newtheorem{theorem}{Theorem}[section]
\newtheorem{proposition}[theorem]{Proposition} 
\newtheorem{corollary}[theorem]{Corollary}
\newtheorem{lemma}[theorem]{Lemma}

\newtheorem*{theorem*}{Theorem}
\newtheorem*{conjecture*}{Conjecture}
\newtheorem*{corollary*}{Corollary}
\newtheorem*{proposition*}{Proposition}

\theoremstyle{definition}
\newtheorem{observation}[theorem]{Observation}
\newtheorem{definition}[theorem]{Definition}
\newtheorem{example}[theorem]{Example}

        {\medskip}

\theoremstyle{remark}

\newcommand{\ol}{\varphi}
\newcommand{\pol}{\Phi}
\def\floor#1{\lfloor {#1} \rfloor}
\def\ceil#1{\lceil {#1} \rceil}

\def\UM#1#2{\bigcup_{#1\in #2}} 
  
\def\VEC#1#2#3{#1_{#2},\ldots,#1_{#3}}
\def\C#1{\left|#1\right|}

\def\nul{\varnothing} 
\def\st{\colon\,}   

\def\esub{\subseteq}  \def\esup{\supseteq}
\def\irep{intersection representation}
\def\orep{overlap representation}
\def\rep{representation}
\def\porep{pure overlap representation}
\def\olap{\leftrightarrow}
\def\disj{\Vert}
 
\def\cF{{\mathcal F}}

\title{Overlap Number of Graphs}
\author{
Daniel W. Cranston\thanks{Virginia Commonwealth University, dcranston@vcu.edu}\,,
Nitish Korula\thanks{University of Illinois, nkorula2@illinois.edu}\,,
Timothy D. LeSaulnier\thanks{University of Illinois, tlesaul2@illinois.edu}\,,
Kevin Milans\thanks{University of Illinois, milans@math.uiuc.edu}\,,\\
Christopher Stocker\thanks{University of Illinois, stocker2@uiuc.edu}\,,
Jennifer Vandenbussche\thanks{Southern Polytechnic State University,
jvandenb@spsu.edu}\,,
Douglas B. West\thanks{University of Illinois, west@math.uiuc.edu.  Work
supported in part by the NSA under award H98230-06-1-0065.}}
\date{\today}

\begin{document}

\maketitle

\begin{abstract}
An {\it overlap representation} of a graph $G$ assigns sets to vertices so that
vertices are adjacent if and only if their assigned sets intersect with neither
containing the other.  The {\it overlap number} $\ol(G)$ (introduced by Rosgen)
is the minimum size of the union of the sets in such a representation.  We
prove the following:
(1) An optimal overlap representation of a tree can be produced in linear time,
and its size is the number of vertices in the largest subtree in which the
neighbor of any leaf has degree 2.
(2) If $\delta(G)\ge 2$ and $G\ne K_3$, then $\ol(G)\le \C{E(G)}-1$, with
equality when $G$ is connected and triangle-free and has no star-cutset.
(3) If $G$ is an $n$-vertex plane graph with $n\ge5$, then $\ol(G)\le 2n-5$,
with equality when every face has length 4 and there is no star-cutset.
(4) If $G$ is an $n$-vertex graph with $n\ge 14$, then
$\ol(G)\le \floor{n^2/4-n/2-1}$, and this is sharp (for even $n$, equality
holds when $G$ arises from $K_{n/2,n/2}$ by deleting a perfect matching).
\end{abstract}

\section{Introduction}

Intersection representations of graphs have been studied for many years. 
An {\it intersection representation} of a graph is a family of sets
corresponding to the vertices so that vertices are adjacent if and only if
their assigned sets intersect.  The first such model was that of {\it interval
graphs}, in which the assigned sets are intervals on the real line.

Intersection representations may use various types of sets.  Erd\H{o}s,
Goodman, and P\'osa~\cite{EGP} introduced intersection representations using
finite sets.  The {\it intersection number} $\theta_1(G)$ is the minimum size
of the union of the sets in an intersection representation of $G$ by finite
sets (\cite{ChW} and~\cite{Eat} use this notation).  In~\cite{EGP}, it was
shown that $\theta_1(G)$ also equals the minimum number of complete
subgraphs needed to cover $E(G)$.

The ``overlap'' model for graph representations arose much later and is less
well studied.  A set {\it overlaps} another set if they intersect but neither
contains the other.  An {\it overlap representation} of a graph $G$ is an
assignment $f$ of sets to the vertices of $G$ so that $uv\in E(G)$ if and only
if $f(u)$ and $f(v)$ overlap.

Just as intersection representations were first studied using intervals, so too
an {\it overlap graph} was defined to be a graph having an overlap
representation using intervals.  The concept appears in the classic book by
Golumbic~\cite{Gol}, noting that a graph is an overlap graph if and only if it
has an intersection representation using chords of a circle.  MathSciNet
returns less than 50 items for ``overlap graph'' and more than 600 for
``interval graph'', though it should be noted that overlap graphs are also
discussed under equivalent terms like ``circle graph''.

Rosgen~\cite{Ros} studied overlap representations using finite sets.  Under any
adjacency rule for assigned sets (such as intersection, containment, or
overlap), a {\it finite representation} of a graph $G$ is a representation in
which the assigned sets are finite.  The {\it size} of a finite representation
$f$ of $G$, denoted $\C f$, is the size of the union of the assigned sets.  The
{\it overlap number} $\ol(G)$ is the minimum size of a finite overlap
representation of $G$.

Throughout this paper, we take $n$ to be the number of vertices of a graph $G$
whose overlap number is being studied.  Rosgen~\cite{Ros} obtained upper bounds
on $\ol(G)$ for trees ($n+1$), chordal graphs ($2n$), planar graphs
($\frac{10}{3}n-6$), and arbitrary graphs ($\theta_1(G)+n$, which yields
$\ol(G)\le \floor{n^2/4}+n$).  He observed that $\ol(K_n)$ is the minimum $t$
such that a $t$-set contains $n$ pairwise incomparable sets, that
$\ol(C_n)=n-1$, and that the overlap number of any caterpillar (with $n>2$) is
the number of vertices in the longest path.  He asked for the maximum value
of $\ol(G)$ in terms of $n$ for trees, chordal graphs, planar graphs, and
arbitrary $n$-vertex graphs, and also for the complexity of computing $\ol$
on trees and on general graphs.

We answer Rosgen's questions about trees using a special subtree.  A {\it
skeleton} is a tree in which the neighbor of any leaf vertex has degree 2.
The largest skeleton in a tree $T$ is unique up to isomorphism, obtained by
deleting all leaves (yielding the {\it derived tree} $T'$) and then restoring
one leaf neighbor of each leaf of $T'$.  Hence we call this {\it the} skeleton
of the tree.  For $n\ge3$, we prove that the overlap number of a tree is the
number of vertices in its skeleton, using an algorithm that produces an \orep\
of this size in linear time.

In Section 2 we give the algorithm and formula for $\ol$ on $n$-vertex trees.
Section 3 presents bounds in terms of the number of edges; we prove that
$\ol(G)\le \C{E(G)}-1$ when $\delta(G)\ge2$ and $G\ne K_3$.  Furthermore, 
equality holds when $G$ is connected, triangle-free, and has no star-cutset,
where a {\it star-cutset} is a separating set $S$ having a vertex $x$ adjacent
to all of $S-\{x\}$.  The results in terms of $\C{E(G)}$ are applied
to $n$-vertex planar graphs in Section 4 and to the family of all $n$-vertex
graphs in Section 5.

In particular, if $G$ is an $n$-vertex plane graph with $n\ge5$, then
$\ol(G)\le 2n-5$, with equality when every face is a 4-cycle and there is no
star-cutset.  When $n\ge14$, the maximum over all $n$-vertex graphs is
$\floor{n^2/4-n/2-1}$, achieved for even $n$ by the graph obtained by deleting
a perfect matching from $K_{n/2,n/2}$.

We note that Henderson~\cite{Hen} independently obtained results on the problems
discussed here.  He obtained constant-factor approximation algorithms for
computing the overlap number on trees and on planar graphs, and he proved that
the maximum overlap number grows quadratically in the number of vertices for a
class of bipartite graphs.  It remains open whether finding the overlap number
is NP-hard in general.

The results in Section 3 use a related model.  A {\it pure \orep} of $G$ is an
\orep\ in which no assigned set contains another.  The {\it pure overlap
number} $\pol(G)$ is the minimum size of a finite pure overlap representation
of $G$.  (Rosgen used the term ``containment-free overlap representation'' for
this model.)  Note that a pure overlap representation of $G$ is both an overlap
representation and an intersection representation of $G$; thus always
$\ol(G)\le\pol(G)$ and $\theta_1(G)\le\pol(G)$.  For this reason, $\pol(G)$
is helpful in proving upper bounds.  Note also that $\rho(H)\le \rho(G)$ when
$\rho\in\{\ol,\pol,\theta_1\}$ and $H$ is an induced subgraph of $G$, since a
representation of $G$ restricts to a representation of $H$.

We say that the vertices adjacent to a vertex $v$ in $G$ are its
{\it neighbors}.  The number of neighbors is the {\it degree} of $v$, denoted
$d_G(v)$ or simply $d(v)$.  The set of neighbors is the {\it neighborhood} of
$v$, denoted $N_G(v)$ or simply $N(v)$.  The {\it closed neighborhood} of $v$,
denoted $N[v]$, is $N(v)\cup\{v\}$.  The minimum vertex degree is $\delta(G)$.
A vertex of degree 1 is a {\it leaf}.  A graph is {\it nontrivial} if it
has at least one edge.

Before beginning the discussion of trees, we prove a lemma used in the lower
bound arguments.  It restricts the form of \orep s.  The idea is due to
Rosgen~\cite{Ros}.

\begin{lemma}\label{containment} 
Let $f$ be an overlap representation of a graph $G$.  If $v\in V(G)$ and
$H$ is a nontrivial component of $G-N[v]$, then either $f(v)$ properly contains
all sets assigned to the vertices of $H$ or $f(v)$ is disjoint from all sets
assigned to vertices of $H$.
\end{lemma}
\begin{proof}
Since no sets used in $H$ overlap $f(v)$, and $H$ is connected, it suffices to
show that if $f(v)\esup f(u)$ for some $u\in V(H)$, and $x\in N(u)$, then
$f(v)\supset f(x)$.

Since $f(u)$ and $f(x)$ overlap, $f(v)\cap f(x)\ne\nul$.  Since $x\notin N(v)$,
we have $f(v)$ and $f(x)$ ordered by inclusion.  Since $x\in N(u)$ forbids
$f(x)\esup f(v)\esup f(u)$, we have $f(v)\supset f(x)$.
\end{proof}

\section{The overlap number of trees}

Rosgen~\cite{Ros} proved that $\ol(T)\le n+1$ when $T$ is a tree.  In fact,
this bound is sharp only for $K_2$.  We provide a linear-time algorithm for
producing an \orep\ of a tree.  We then prove that this \rep\ is optimal.

A {\it caterpillar} is a tree in which all edges are incident to a single path.
Rosgen~\cite{Ros} proved that the overlap number of any caterpillar equals the
number of vertices in a longest path.  For a caterpillar, this path is the 
skeleton.  We will need this result along with a technical property of the
representation, because our procedure for extending a \rep\ along an added
caterpillar differs from the \rep\ for the initial caterpillar.  

\begin{definition}\label{minimal}
For an \orep\ $f$ of a graph $G$, the {\it associated poset} $P_f$ is the
inclusion order on $\{f(v)\st v\in V(G)\}$.  A vertex $v$ is {\it minimal} in
$f$ if $f(v)$ is a minimal element of $P_f$, and $v$ is {\it $a$-minimal} if
$f(v)$ is a minimal element of the subposet of $P_f$ consisting of the elements
that contain $a$.
In the same way that $\esup$ means ``contains'', we use $\olap$ to mean
``overlaps'' and ``$\disj$'' to mean ``does not intersect''.
\end{definition}

\begin{lemma}\label{cater}
  Let $T$ be a caterpillar whose longest path has vertices $\VEC v1l$
  in order.  If $l\ge 3$, then $T$ has an \orep\ $f$ of size $l$.
  Furthermore, with $\{\VEC a1l\}$ being the union of the assigned sets,
  $f$ may be chosen so that $v_i$ is $a_i$-minimal in $1\le i\le l-1$.
\end{lemma}
\begin{proof}
Let $f(v_i)=\{a_i,a_{i+1}\}$ for $1\le i\le l-1$.  All leaves (including
$v_l$) have a neighbor in $\{\VEC v2{l-1}\}$.  For each leaf neighbor $x$
of $v_i$, let $f(x)=\{\VEC a1i\}$.

By construction, $f(v_{i-1})\olap f(v_i)$ for $2\le i\le l-1$, and
nonconsecutive sets in that list are disjoint.  If $x$ is a leaf neighbor of
$v_i$, then $f(x)\olap f(v_i)$, $f(x)\esup f(v_j)$ for $j<i$, and
$f(x)\disj f(v_j)$ for $j>i$.  Also the sets assigned to leaves form a chain
by inclusion.  Hence $f$ is an \orep\ of $T$.  Since no assigned sets are
singletons, $v_i$ is $a_i$-minimal.
\end{proof}

\begin{observation}\label{addelt}
If $A$ and $B$ are sets such that $A\esup B$ or $A\olap B$, then adding
an element not in $A\cup B$ to $A$ or to both $A$ and $B$ preserves the
relation.  If $A\disj B$, then the relation is preserved when the 
element is added to just one of $\{A,B\}$.
\end{observation}

\begin{lemma}\label{catadd}
  Let $G$ be the union of a graph $H$ and a caterpillar $T$ such that $H\cap T$
  consists of one vertex $v$ that is not isolated in $H$ and is an endpoint of a
  longest path in $T$.  Let $H$ have an \orep\ $f$ of size $m$, and
  let $\VEC w0l$ be the vertices along a longest path in $T$, with $v=w_0$.
  If $v$ is $a$-minimal in $f$ for some $a\in f(v)$, then $G$ has an \orep\ $f'$
  of size $m+l$, with added elements $\VEC b1l$, such that $w_i$ is
  $b_i$-minimal in $f'$ for $1\le i\le l$, and any vertex of $H$ that
  is $c$-minimal in $f$ is also $c$-minimal in $f'$.
\end{lemma}
\begin{proof}
Let $B=\{\VEC b1l\}$, and let $b_0=a$.  Let $f'(v)=f(v)$, and let
$f'(w_i)=\{b_{i-1},b_i\}$ for $1\le i\le l-1$.  Each remaining vertex of $T$ is
a leaf with neighbor in $\{\VEC w1{l-1}\}$.  For each leaf $x$ in $T$ with
neighbor $w_i$, let $f'(x)=\{\VEC bil\}$.  For $u\in V(H)-\{v\}$,
let $f'(u)=f(u)$ if $a\notin f(u)$; otherwise, let $f'(u)=f(u)\cup B$.

By construction, $f'$ generates a path on $\VEC w0l$, since $d_H(v)\ge1$
requires $f(v)\ne\{a\}$.  If $x$ is a leaf in $T$ adjacent to $w_j$, then
$f'(x)$ contains the sets assigned to $w_i$ and its leaf neighbors for $i>j$.
Also $f'(x)\olap f'(w_j)$, and $f'(x)\disj f'(w_i)$ for $i<j$.

If $u\in V(H)-\{v\}$ and $y\in V(T)-\{v\}$, either $f'(u)\disj f'(y)$ or
$f'(u)\esup f'(y)$, depending on whether $f'(u)$ acquires $B$.  Since
$B\esub f'(u)$ if and only if $a\in f(u)$, by Observation~\ref{addelt} the
relation between sets assigned to vertices of $V(H)-\{v\}$ under $f'$ and $f$
is the same.

Since $f'(v)\disj B$, among the sets assigned by $f'$ to $V(T)-\{v\}$ only
$f'(w_1)$ overlaps $f'(v)$.  Now compare $v$ with a vertex $u\in V(H)-\{v\}$.
Since $v$ is $a$-minimal, $f(u)\subset f(v)=f'(v)$ implies $f'(u)=f(u)$.
Otherwise, Observation~\ref{addelt} implies that $f'(u)$ and $f'(v)$ have the
same relation as $f(u)$ and $f(v)$.  We have shown that $f'$ is an overlap
representation of $G$.

Note that as in Lemma~\ref{cater}, each $w_i$ is $b_i$-minimal in $f'$.  If $u$
is $c$-minimal in $f$, then for every $y$ with $c\in f(y)$,
Observation~\ref{addelt} implies in all cases that $u$ is $c$-minimal in $f'$.
\end{proof}


\begin{theorem}\label{treeupper}
Every tree other than $K_2$ has an \orep\ whose size is the number of
vertices in its skeleton.
\end{theorem}
\begin{proof}
We grow a tree $T$ by successive addition of appropriate caterpillars.  The
first caterpillar, $T_0$, consists of any maximal subtree of $T$ that is a
caterpillar whose leaves are also leaves of $T$.  The maximality guarantees
that the ends of a longest path in $T_0$ are leaves of $T$ that are also leaves
in the skeleton.

When the subtree absorbed so far is $T_i$, the next caterpillar $T'$ is a
maximal caterpillar contained in $T$ such that an endpoint $x$ of some longest
path of $T'$ (and no other vertex of $T'$) is in $T_i$, and all leaves of $T'$
are leaves in $T$.  Let $T_{i+1}= T_i\cup T'$.  The end opposite $x$ of a
longest path in $T'$ is a leaf of $T$ that is preserved in the skeleton.  Thus
the maximality conditions guarantee that the subtree of $T$ formed by the union
of the longest paths in the chosen caterpillars is the skeleton of $T$.

By Lemma~\ref{cater}, the initial caterpillar has an \orep\ of the desired
size, with all non-leaf vertices being $c$-minimal for distinct choices of $c$.
By Lemma~\ref{catadd}, the process continues with the $b$-minimality conditions
on non-leaf vertices preserved and the desired number of elements being added
at each step.  (In fact, in the final \orep\ $f$, only one vertex of the 
skeleton is not $c$-minimal for any $c$; it is a leaf of $T_0$.) 
\end{proof}

The skeleton of any tree $G$ is an induced subgraph of $G$.  Therefore, to
prove that the \orep\ produced in Theorem~\ref{treeupper} is optimal for every
tree with $n\ge3$, it suffices to show that if $T$ is a skeleton with $n$
vertices, then $\ol(T)=n$.

The idea of the proof is inductive.  Given an \orep\ $f$ for a skeleton $T$, we
seek one or two vertices in $T$ (a leaf or a leaf and its neighbor) whose
deletion yields a smaller skeleton $T'$ for which we can obtain an \orep\ by
deleting one or two elements from $f$.  The lower bound then follows
inductively.  To do this, we need to know when elements can be deleted from an
\orep\ $f$ or from a restriction of $f$ to a subgraph.  We write $f-S$ for the
result of subtracting $S$ from each set assigned under $f$.

\begin{lemma}\label{deletion}
If $f$ is an \orep\ of a graph $G$, then $f-S$ is an \orep\ of $G$ if and only
if $S$ does not contain the intersection or difference of the sets assigned
to any two adjacent vertices of $G$.
\end{lemma}

\begin{proof}
{\it Necessity:}
Deleting a set containing the intersection or difference of the sets for
adjacent vertices would delete that edge from the corresponding overlap graph.

{\it Sufficiency:}
Deleting a set $S$ satisfying the stated condition maintains the overlap
condition for any pair of overlapping sets.  Deletions from disjoint sets
maintain disjointness, and containments are preserved because $A\esub B$
implies $A-S\esub B-S$.
\end{proof}

\begin{definition}
Let $f$ be an assignment of sets to $V(G)$.  A set $S$ of elements is
{\it $f$-uniform} if every assigned set $f(v)$ contains all or none of $S$.
\end{definition}

\begin{observation}\label{uniform}
If $f$ is an \orep\ of a graph $G$, then every proper subset of an $f$-uniform
set is deletable from $f$.  Hence an \orep\ having a uniform set of size 2 is
not optimal.
\end{observation}

Our next lemma is the key tool in proving the lower bound for trees.
It strengthens Observation~\ref{uniform}, allowing us to reduce the size of an
overlap representation when it has a set that is uniform except at one vertex.

\begin{lemma}\label{hammer-lem}
Let $v$ be a vertex in a graph $G$ such that $N(v)$ is independent and contains
no leaves.  Let $f$ be an overlap representation of $G$, and let $f'$ be its
restriction to $G-v$.  If $\{a,b\}$ is $f'$-uniform, then $f-\{a\}$ or
$f-\{b\}$ is an overlap representation of $G$.
\end{lemma}

\begin{proof}
By Observation~\ref{uniform}, the claim follows unless exactly one element of
$\{a,b\}$ is in $f(v)$.  Hence we may assume that $a\notin f(v)$ and
$b\in f(v)$.

Suppose that $f-\{a\}$ is not an overlap representation of $G$.
Since Observation~\ref{uniform} implies that $f-\{a\}$ is an overlap \rep\ of
$G-v$, Lemma~\ref{deletion} implies that some edge incident to $v$ is lost
when $a$ is deleted from $f$.  Let $vw_1$ be such an edge.  Because
$a\notin f(v)$, we conclude that $f(w_1)-f(v)=\{a\}$.  Because $\{a,b\}$ is 
$f'$-uniform, also $b\in f(w_1)$.

If $f-\{b\}$ also is not a representation, then deleting $b$ also destroys some
edge $vw_2$ incident to $v$.  Since $b\in f(v)$, either $f(v)\cap f(w_2)=\{b\}$
or $f(v)-f(w_2)=\{b\}$.  We obtain a contradiction from each case.
Note first that since each $w_i$ has a neighbor other than $v$, and $\{a,b\}$
is $f'$-uniform, each $f(w_i)$ contains an element outside $\{a,b\}$.

{\it Case 1: $f(v)\cap f(w_2)=\{b\}$.}  Since $\{a,b\}$ is $f'$-uniform,
also $a\in f(w_2)$.  Thus $\{a,b\}\esub f(w_1)\cap f(w_2)$.
Since $w_1w_2\notin E(G)$, the sets $f(w_1)$ and $f(w_2)$ cannot each have an
element outside the other.  However, each has an element outside $\{a,b\}$, so
they share another element $c$.  Now $f(w_1)-f(v)=\{a\}$ yields $c\in f(v)$,
while $f(v)\cap f(w_2)=\{b\}$ yields $c\notin f(v)$.

{\it Case 2: $f(v)-f(w_2)=\{b\}$.}
Since $\{a,b\}$ is $f'$-uniform, $f(w_2)\cap\{a,b\}=\nul$.  Since
$f(w_1)-f(v)=\{a\}$, the guaranteed element $c$ in $f(w_1)-\{a,b\}$ must lie in
$f(v)$.  Since $f(v)-f(w_2)=\{b\}$, also $c\in f(w_2)$.  Meanwhile,
$f(v)\olap f(w_2)$ requires an element $d$ in $f(w_2)-f(v)$.  Since
$f(w_1)-f(v)=\{a\}$, we have $d\notin f(w_1)$.  Now $f(w_1)$ and $f(w_2)$
share $c$ and overlap, contradicting $w_1w_2\notin E(G)$.
\end{proof}

In a skeleton, the neighbor of a leaf vertex has degree 2.

\begin{definition}
In an \orep\ $f$ of a skeleton $T$, a leaf $l$ is {\em doubly-minimal} if both
$l$ and the neighbor of $l$ are minimal in $f$.
\end{definition}

\begin{lemma}\label{doubly-minimal-leaf}
In an \orep\ $f$ of any skeleton $T$, there is at most one nonminimal leaf.
If $T\ne P_4$, then there is a doubly-minimal leaf.
\end{lemma}
\begin{proof}
If $T=P_3$, then there is only one non-edge, so only for the two leaves can
one set properly contain another.  Thus a leaf and the center are minimal in
$f$.  Henceforth assume $T\ne P_3$.  In a skeleton other than $P_3$, no two
neighborhoods are equal.  Hence also no two assigned sets are equal.

Since the neighbor of any leaf $x$ has degree 2, $T-N[x]$ is connected.  If $x$
is nonminimal, then by Lemma~\ref{containment} $f(x)$ properly contains the
sets assigned to all vertices other than its neighbor, including the other
leaves.  Hence only one leaf can be nonminimal.

Let $A$ be the set of neighbors of minimal leaves; we have shown that
$A\ne\nul$.  Choose $v\in A$ such that $f(v)$ is minimal in 
$\{f(y)\st y\in A\}$.  Let $N(v)=\{x,u\}$, with $x$ being the leaf.
Since $T\ne P_3$, $u$ is not a leaf.

When $T=P_4$, the claim fails when the sets in $f$ are $ab$, $bceg$, $abde$,
and
$eg$.

If $T\ne P_4$, then $u$ has no leaf neighbor, and each component of $T-N[v]$ is
nontrivial.  If $x$ is not doubly-minimal, then $f(v)$ properly contains the
sets for all vertices in some component $T'$ of $T-N[v]$, by
Lemma~\ref{containment}.  Let $x'$ be a leaf of $T$ contained in $T'$, and let
$v'$ be its neighbor, also in $T'$.  Since $f(v)$ contains $f(v')$, the choice
of $v$ from $A$ requires $x'$ to be nonminimal.  As observed earlier, this
yields $f(v)\subset f(x')$, contradicting $f(v)\supset f(x')$.
\end{proof}

\begin{theorem}\label{lowerbound}
If $T$ is a skeleton with $n$ vertices, where $n\ge3$, then $\ol(T) \geq n$.
\end{theorem}
\begin{proof}
We first note a tool that allows us to apply Lemma~\ref{hammer-lem} when a
leaf is a minimal vertex.  Recall that $N[v]$ denotes $N(v)\cup\{v\}$.

\smallskip
{\narrower
\it
\noindent
(*) If $x$ is a minimal vertex in an \orep\ $f$ of a graph $G$, and $f'$ is the
restriction of $f$ to $G-N[x]$, then $f(x)$ is $f'$-uniform.

}

\smallskip
\noindent
To prove (*), note that $v\in V(G)-N[x]$ implies $xv\notin E(G)$.  Hence
$f(x)\disj f(v)$ or $f(x)\esub f(v)$ or $f(x)\supset f(v)$.  Minimality of
$x$ excludes the last, so $f(v)$ contains all or none of $f(x)$.

We prove the lower bound on $\ol(T)$ by induction on $n$.  Since $P_3$ has an
edge, $\ol(P_3)\ge3$, so we may assume $n\ge4$.  Let $f$ be an optimal overlap
representation of $T$; Lemma~\ref{doubly-minimal-leaf} yields a leaf $x$ of $T$
that is minimal in $f$.  Let $v$ be the neighbor of $x$, and let $u$ be the
other neighbor of $v$; note that $d(u)\ge2$.  Since $xv \in E(G)$, there exist
$a\in f(x)-f(v)$, $b\in f(x)\cap f(v)$, and $c \in f(v)-f(x)$.

Let $T'=T-x$ and $T''=T-x-v$.  Let $f'$ and $f''$ be the restrictions of $f$ to
$T'$ and $T''$, respectively.  We consider two cases, depending upon $d(u)$.

If $d(u)=2$, then $T'$ is a skeleton.  Since $x$ is minimal, (*) implies that
$f(x)$ (and therefore $\{a,b\}$) is $f''$-uniform.  Since $d(u)=2$, the
neighborhood of $v$ in $T'$ is independent and contains no leaves.  Hence
Lemma~\ref{hammer-lem} applies, and $f'-\{a\}$ or $f'-\{b\}$ is an overlap
representation of $T'$.  By the induction hypothesis, $|f'|\ge n-1$, so
$|f|\ge n$.

If $d(u)>2$, then $T\ne P_4$.  Now Lemma~\ref{doubly-minimal-leaf} allows us to
choose $x$ to be doubly-minimal in $f$.  Since $d(u)\ge3$, deleting $x$ and $v$
from $T$ does not create any new leaves, so $T''$ is a skeleton.  Since $x$ is
minimal, (*) implies that $f(x)$ (and therefore $\{a,b\}$) is $f''$-uniform.
Thus $f''-\{a\}$ is an overlap representation of $T''$.  Let $g=f''-\{a\}$.

Since $x$ is doubly-minimal, $v$ is minimal in $f$, and thus (*) implies that
$f(v)$ (and therefore $\{b,c\}$) is $g$-uniform.  We now apply
Lemma~\ref{hammer-lem} to the vertex $u$, graph $T''$, and \orep\ $g$ of $T''$.
Let $g'$ be the restriction of $g$ to $T''-u$.  Since $\{b,c\}$ is
$g'$-uniform, and $d(u) \geq 3$ implies that $u$ has no leaf neighbors in the
skeleton $T''$, Lemma~\ref{hammer-lem} implies that $g-\{b\}$ or $g-\{c\}$ is
an overlap representation of $T''$.  By the induction hypothesis, $\C g\ge n-1$
and $\C f\ge n$.
\end{proof}

We have proved the following conclusion.

\begin{theorem}\label{tree}
If $T$ is a tree, then $\ol(T)$ is the number of vertices in the skeleton of
$T$.  Furthermore, there is a linear-time algorithm to produce an optimal
\orep.
\end{theorem}

\section{Bounds from the Number of Edges}

As mentioned earlier, Erd\H{o}s, Goodman, and P\'osa~\cite{EGP} observed that
finite intersection \rep s of a graph $G$ correspond to families of complete
subgraphs covering $E(G)$.  In cases where the \irep\ arising from a
decomposition into complete subgraphs is also an \orep, its size must be at
least the pure overlap number $\pol(G)$, defined in Section 1.  On the other
hand, always $\ol(G)\le\pol(G)$.

Several upper bounds will be used repeatedly in the remainder of the paper;
we give them names to improve readability.  A {\it decomposition} of a graph
$G$ is a family $\cF$ of pairwise edge-disjoint subgraphs whose union is $G$.

\begin{lemma}[\rm Decomposition Bound]\label{clique-decomp}
  Let $\cF$ be a decomposition of a graph $G$ into complete subgraphs of order
  at most $k$, where $k\ge2$.  If $\delta(G)\ge k$, then $\pol(G) \leq |\cF|$.
  In particular, $\delta(G)\ge2$ implies $\pol(G)\le \C{E(G)}$.
\end{lemma}
\begin{proof}
For each $v\in V(G)$, let $f(v)$ be the set of all members of $\cF$ that
contain $v$.  Each edge lies in some member of $\cF$, so $f$ is an \irep.

A vertex has at most $k-1$ neighbors in a complete subgraph of order $k$.
Since $\delta(G)\ge k$, each $|f(v)|$ is at least $2$.  Each edge is covered
only once, so $\C{f(v)\cap f(u)}\le 1$.  Hence no assigned set contains
another, and $f$ is a pure \orep.
\end{proof} 

\begin{corollary}\label{edgebound} 
If $G$ is triangle-free, then $\pol(G)\ge\C{E(G)}$, with equality when
$\delta(G)\ge 2$.
\end{corollary}
\begin{proof}
A \porep\ is also an \irep.
\end{proof}

Lemma~\ref{clique-decomp} provides an upper bound when $\delta(G)\ge2$, and we
will also apply it with $k=3$ for decompositions into edges and triangles.
Hence we want vertices of degree at most 2 to contribute little to $\pol(G)$.

\begin{lemma}[\rm Deletion Bound]\label{poldeg2}
  If $v$ is a vertex of degree at most 2 in a graph $G$ with at least
  three vertices, then $\pol(G)\le \pol(G-v)+2$, with strict
  inequality when $d(v)=0$.
\end{lemma}
\begin{proof}
If $d(v)=0$, then to avoid overlap and containment with all other assigned
sets, we must assign $v$ an element not assigned to any other vertex.  Thus
$\pol(G)=\pol(G-v)+1$.

For $d(v)\in\{1,2\}$, let $f$ be an optimal \porep\ of $G-v$.  Introduce new
labels $a$ and $b$.  Let $f'(v)=\{a,b\}$.  Let $f'(x)=f(x)$ for $x\notin N[v]$.
For $x\in N(v)$, let $f'(x)=f(x)\cup \{c\}$, where $c$ is one of the new
labels, each used once when $d(v)=2$.

Changing from $f$ to $f'$ creates no new intersections except to establish the
edge(s) incident to $v$.  Adding a new element to its neighbor(s) does not
create containments, and there is no containment involving $f(v)$ and a set
assigned to a neighbor, since the neighbors also receive an old label (even if
isolated in $G-v$).
\end{proof}

Next we discuss bounds on $\ol$ in terms of the number of edges.  In
contrast to $\pol(G)$, generally $\ol(G)<|E(G)|$ (though not for skeletons, as
we have seen).  An easy reduction allows us to forbid repeated vertex
neighborhoods and isolated vertices.

\begin{observation}\label{freevertices}
If a graph $G$ has a vertex $v$ such that $N(v)$ is empty or equals another
vertex neighborhood, then $\ol(G)=\ol(G-v)$.  In the first case, we extend an
overlap representation $f$ of $G-v$ by assigning $v$ the set
$\UM u{V(G-v)} f(u)$.  In the second case, we assign $v$ the same set as $w$,
where $N(v)=N(w)$.
\end{observation}

Let $B_n$ denote the graph that is the union of $n-2$ triangles having
a common edge; this graph is sometimes called the {\it $n$-book}.

\begin{lemma}\label{triangles}
  The overlap number of the $n$-book $B_n$ is $3$.
\end{lemma}
\begin{proof}
Since all the vertices besides the two vertices of degree $n-1$ have identical
neighborhoods, Observation~\ref{freevertices} allows us to remove them without
changing the overlap number until we are left with a triangle, which has
overlap number 3.
\end{proof}

\begin{lemma}[\rm Edge Bound]\label{ol-edge-upper}
  Let $G$ be an $n$-vertex graph other than the book $B_n$.  If $\delta(G)\ge2$
  and $uv\in E(G)$, then $G$ has an overlap representation $f$ with size
  $\C{E(G)}-1$ such that neither $f(u)$ nor $f(v)$ is properly contained in the
  set assigned to any other vertex.  In particular, $\ol(G)\le \C{E(G)}-1$ when
  $\delta(G)\ge2$ unless $G$ is the $3$-book $K_3$. 
\end{lemma}
\begin{proof}
We define an explicit representation using a label for each edge {\it other
than} $uv$.  For $w\notin\{u,v\}$, let $f(w)$ be the set of labels for edges
incident to $w$.  For $w\in\{u,v\}$, let $f(w)$ be the set of labels for edges
{\it not} incident to $w$.  The restriction of $f$ to $G-u-v$ is a \porep\ of
$G-u-v$ (labels for edges to $u$ or $v$ can establish non-containment).

By construction, $f(u)\esup f(w)$ when $w$ is a nonneighbor of $u$.  This
establishes nonadjacency and prohibits $f(u)$ from proper containment in
another assigned set.  Similarly for $v$.  (However, $f(u)=f(w)$ when
$G=K_{2,n-2}$ and $\{u,w\}$ is a partite set of size $2$.)

For $w\in N(u)-\{v\}$, the label for $uw$ is in $f(w)-f(u)$.  Since $d(w)\ge2$,
the label for some other edge incident to $w$ lies in $f(u)\cap f(w)$.
To establish $f(u)-f(w)\ne\nul$, it suffices to have an edge incident to
neither $w$ nor $u$.  If every edge is incident to $w$ or $u$, then $G=B_n$,
and $\ol(G)=3$.  The same argument applies to edges at $v$.

For the edge $uv$ itself, $f(u)-f(v)$ contains the label for an edge other than
$uv$ incident to $v$.  Similarly, $f(v)-f(u)\ne \nul$.  To ensure that
$f(u)\cap f(v)\ne\nul$, we need an edge incident to neither $u$ nor $v$.  As
above, this exists unless $G=B_n$.
\end{proof}

In Theorem~\ref{ol-edge-lower}, we will prove equality in the upper bound
$\ol(G)\le \C{E(G)}-1$ for a special family of graphs with $\delta(G)\ge2$.
For this we will need a definition and two lemmas.

\begin{definition}
  A {\it star-cutset} in a graph $G$ is a separating set $S$ containing a
  vertex $x$ adjacent to all of $S-x$.  If $G$ has no star-cutset, then it is
  {\em star-cutset-free}.
\end{definition}

\begin{lemma}\label{star-cutset}
  If $f$ is an \orep\ of a connected graph $G$ with no star-cutset,
  then any two vertices that are not minimal in $f$ are adjacent.
\end{lemma}
\begin{proof}
Let $u$ and $v$ be such vertices.  If $v\notin N(u)$, then $v$ remains in
$G-N[u]$.  Since $G$ has no star-cutset, $G-N[u]$ is connected.  Since $u$ is
not minimal, $f(u)$ properly contains the set assigned to some vertex of
$G-N[u]$.  By Lemma~\ref{containment}, $f(u)$ properly contains $f(v)$.
The same argument yields $f(v)\supset f(u)$, a contradiction.
\end{proof}

\begin{lemma}\label{star-cutset-repeated-nbrhd}
  If $G$ is an $n$-vertex triangle-free graph with no star-cutset, then $G$
  does not have distinct vertices with the same neighborhood, unless
  $G=K_{2,n-2}$ with $n\le 4$.
\end{lemma}
\begin{proof}
If $N(u) = N(v)$, then $v$ is isolated in $G-N[u]$.  Since $G$ has no
star-cutset, $G-N[u]$ contains only $v$.  Thus $V(G)=\{u,v\}\cup N(u)$.
Also $N(u)$ induces no edges, since $G$ is triangle-free.  Thus
$G=K_{2,n-2}$.  Also, $n-2\le 2$, since otherwise deleting $N[x]$ for some
$x\in N(u)$ disconnects $G$, contradicting the absence of star-cutsets.
\end{proof}

\begin{theorem}\label{ol-edge-lower}
  If $G$ is a triangle-free graph with no star-cutset, then
  $\ol(G)\ge\C{E(G)}-1$.
\end{theorem}
\begin{proof}
Let $f$ be an \orep\ of $G$.  If two assigned sets are equal, then those
vertices have the same neighborhood, and Lemma~\ref{star-cutset-repeated-nbrhd}
yields $G=K_{2,n-2}$ with $n\le4$.  By Observation~\ref{freevertices},
$\ol(K_{2,0}) = 1$ and $\ol(K_{2,1}) = \ol(K_{2,2}) = 3$.  In each case,
$\ol(G)\ge\C{E(G)}-1$.

Hence we may assume that no two sets assigned by $f$ are equal.  Since $G$ is
triangle-free, by Lemma~\ref{star-cutset} at most two vertices are non-minimal.
We consider three cases.

{\it Case 0: Every vertex is minimal in $f$.}  In this case, $f$ is a \porep,
and Corollary~\ref{edgebound} yields $\C f\ge \C{E(G)}$.

{\it Case 1: One vertex $u$ is nonminimal in $f$.}  Since $G-N[u]$ is
connected and $f(u)$ contains some other assigned set, $f(u)$ contains all
elements assigned to the nonneighbors of $u$.  Also $N(u)$ is independent,
so every edge of $G-u$ has an endpoint outside $N[u]$.

All containments involve $f(u)$.  Hence $f$ restricts to a \porep\ and thus an
\irep\ on $G-u$.  Since $G$ is triangle-free, for each edge $e$ of $G-u$ there
is an element assigned by $f$ to the endpoints of $e$.  It also lies in $f(u)$,
since $e$ has an endpoint outside $N[u]$.  Let $S$ be this set of elements.

Since $S\esub f(u)$, we still must make $f(w)-f(u)$ nonempty for $w\in N(u)$.
Since $N(u)$ is independent and $u$ is the only nonminimal vertex, the sets
assigned to $N(u)$ are pairwise disjoint.  Hence $N(u)$ requires distinct
additional elements, yielding $\C f\ge\C{E(G)}$.

{\it Case 2: Two vertices, $u$ and $v$, are nonminimal in $f$.}  By
Lemma~\ref{star-cutset}, $uv\in E(G)$.  As above, $f$ restricts to an \irep\ on
$G-u-v$ with an element for each edge; let $S$ be this set of elements.
Since $u$ is nonminimal and $G-N[u]$ is connected, $f(u)$ contains all elements
assigned to vertices outside $N[u]$.

As above, each $w\in N(u)-\{v\}$ needs an element not in $f(u)$, and these
elements are distinct since $G$ is triangle-free.  The same holds for $N(v)$.
We thus have $\C f\ge\C{E(G)}-1$ unless there exist $x\in N(u)-\{v\}$ and
$y\in N(v)-\{u\}$ with $f(x)$ and $f(y)$ having a common element outside
$f(u)\cup f(v)$.  Since $G$ is triangle-free, $u$ and $v$ have no common
neighbors, so $f(u)\supset f(y)$ and $f(v)\supset f(x)$.  Hence the elements
establishing the edges between $\{u,v\}$ and their neighbors are distinct, and
$\C f\ge\C{E(G)}-1$.
\end{proof}

The proof of Theorem~\ref{ol-edge-lower} shows that for a connected
triangle-free graph $G$ with $\delta(G)\ge2$ and no star-cutset, the {\it only}
way to form an \orep\ with fewer than $\C{E(G)}$ elements is as in the proof of
the Edge Bound (Lemma~\ref{ol-edge-upper}).

Theorem~\ref{ol-edge-lower} determines the overlap number for graphs that we
will show are extremal in the classes of $n$-vertex graphs and $n$-vertex
planar graphs.  In each example below, the graphs are connected and
triangle-free without star-cutsets, and the number of edges is one more than
the specified overlap number (prohibition of star-cutsets requires
$\delta(G)\ge2$).

\begin{corollary}\label{bicliq-match}
  For $n\ge6$, the $n$-vertex graph obtained by deleting from
  $K_{\floor{n/2},\ceil{n/2}}$ a matching of size $\floor{n/2}$ (and deleting
  one edge incident to the remaining vertex when $n$ is odd) has overlap number
  $\floor{n^2/4-n/2-1}$.
\qed
\end{corollary}

\begin{corollary}\label{planar-max}
  If $G$ is a triangle-free plane graph in which every face has length $4$, and
  $G$ has no star-cutset, then $\ol(G)=2n-5$.
\qed
\end{corollary}

\begin{example}\label{planar-examp}
Graphs as described in Corollary~\ref{planar-max} exist for $n\ge10$ (also for
$n=4$ and $n=8$).  When $n\equiv 0\mod 4$, the cartesian product of $P_{n/4}$
and $C_4$ suffices.  When $n\equiv 0\mod 2$, one can start with an even cycle
$C$ and add a vertex inside adjacent to the even-indexed vertices on $C$ and a
vertex outside adjacent to the odd-indexed vertices on $C$.

For odd $n$, take such a graph $G$ with $n-1$ vertices embedded in the plane,
let $x$ be a vertex of degree at least 4 in $G$ ($x$ exists if $\C{V(G)}\ge10$),
and let $u$ and $v$ be nonconsecutive neighbors of $x$ in the embedding.  Form
$G'$ by replacing $x$ with nonadjacent vertices $x'$ and $x''$ whose
neighborhoods in the new graph $G'$ partition $N_G(x)$, except that
$N_{G'}(x')\cap N_{G'}(x'')=\{u,v\}$.  The vertices $\{x',u,x'',v\}$ form a new
face surrounding the former edges $xu$ and $xv$, and the other edges at $x$
attach instead to $x'$ and $x''$.
\qed
\end{example}

As we did with pure overlap number, for overlap number we will want to 
accommodate vertices of degree less than 2 without much cost.  By
Observation~\ref{freevertices}, we may assume that there are no isolated
vertices and that vertex neighborhoods are distinct.

The corresponding result for vertices of degree 1 is a special case of a more
general result (proved in the same way) that permits saving labels when \orep s
of subgraphs are combined at a cut-vertex.  For clarity, we present only the
result that we use to obtain our extremal results in the subsequent sections.

\begin{lemma}\label{ol-leaves} 
If $v$ is a leaf in a graph $G$ and $G-v$ is nontrivial, then
$\ol(G) \le \ol(G-v) + 2$.
\end{lemma}
\begin{proof} 
Let $u$ be the neighbor of $v$.  If $uv$ is an isolated edge and $G-v$ is
nontrivial, then let $R=\UM x{V(G-\{u,v\})} f(x)$, where $f$ is an \orep\ of
$G-\{u,v\}$.  Note that $\C R\ge1$.  Modify $f$ by assigning $R\cup\{a\}$ to
$u$ and $R\cup\{b\}$ to $v$, where $a,b\notin R$.  This produces an \orep\ of
$G$, so $\ol(G)\le \ol(G-\{u,v\})+2\le \ol(G-v)+2$.

Hence we may assume that $d_G(u)\ge2$.  Let $f$ be an optimal \orep\ of $G-v$.
Let $W=V(G)-\{u,v\}$.  Define $f'$ on $V(G)$ as follows.  Let $f'(v)=S=\{a,b\}$,
where $a,b\notin \UM xW f(x)$.  Let $f'(u)=f(u)\cup \{b\}$.  For $x\in W$, let
$f'(x)=f(x)\cup S$ if $f(x)\esup f(u)$; otherwise, let $f'(x)=f(x)$.  Note that
$\C {f'}=\C{f}+2$.

We check that $f'$ is an \orep\ of $G$.  Since $\C{f(u)}\ge 2$, we have
$f'(u)\olap f'(v)$.  For $x\in W$, we have $f'(x)\disj f'(v)$ or
$f'(x)\esup f'(v)$, depending on whether $f'(x)$ acquires $S$, so $v$ receives
no other edges.

For $x,y\in W$, the assigned sets acquire $S$ if and only if they contain
$f(u)$.  By Observation~\ref{addelt}, the relation between $f'(x)$ and $f'(y)$
is the same as between $f(x)$ and $f(y)$.  If $f'(x)=f(x)$, then
$f(x)\not\esup f(u)$, and the relation between $f(x)$ and $f(u)$ is preserved.
If $f(x)\esup f(u)$, then $f'(x)=f(x)\cup S$ and again the relation is
preserved.
\end{proof}

\section{Overlap Number of Planar Graphs}

In order to apply the Decomposition Bound    
for planar graphs that may contain triangles, we need an efficient
decomposition into small complete subgraphs.  By Euler's Formula, a
triangle-free planar graph has as most $2n-4$ edges, with equality only if
every face is a 4-cycle. 

\begin{lemma}\label{plan-decomp}
  If $G$ is an $n$-vertex plane graph, and $n\ge3$, then $G$
  decomposes into at most $2n-5$ edges and facial triangles unless:
\\
\quad
  (a) every face is a $4$-cycle, in which case $G$ decomposes into $2n-4$ edges, or
\\
\quad
  (b) $G$ is $K_4$, which decomposes into three edges and one facial triangle.
\end{lemma}
\begin{proof}
We use induction on the number of facial triangles in $G$.  If there are none,
then Euler's Formula suffices.  If $G$ has a facial triangle $T$, then form a
plane graph $G'$ from $G$ by deleting $E(T)$ and introducing a new vertex $v$
adjacent to $V(T)$.  Since $v$ belongs to no triangle, $G'$ has fewer facial
triangles than $G$.

Suppose first that $G'$ has a facial cycle that is not a 4-cycle.  By the
induction hypothesis, $G'$ decomposes into at most $2n-3$ triangles and edges
($G'$ has $n+1$ vertices).  Since $v$ is in no triangle, the three edges
incident to $v$ are edges in the decomposition.  Replacing them with $T$ yields
the desired decomposition of $G$.

If every face in $G'$ is a 4-cycle, then each edge of $T$ lies in another
facial triangle in $G$.  By the induction hypothesis, $G'$ decomposes into
$2n-2$ edges.  If the faces incident to $v$ in $G'$ have no shared edges
not incident to $v$, then their nine edges $G'$ can be replaced with three
triangles to decompose $G$ into $2n-8$ edges and facial triangles
(Figure~\ref{figplan}a).

If two of these faces share an edge, then the eight distinct edges can be
replaced with two triangles and two edges to decompose $G$ into $2n-6$ edges
and facial triangles (Figure~\ref{figplan}b).

\begin{figure}[hbt]
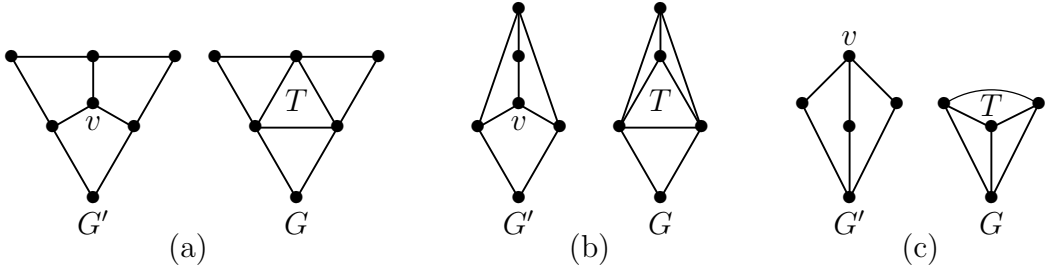

\gpic{
\expandafter\ifx\csname graph\endcsname\relax \csname newbox\endcsname\graph\fi
\expandafter\ifx\csname graphtemp\endcsname\relax \csname newdimen\endcsname\graphtemp\fi
\setbox\graph=\vtop{\vskip 0pt\hbox{%
    \graphtemp=.5ex\advance\graphtemp by 0.309in
    \rlap{\kern 0.490in\lower\graphtemp\hbox to 0pt{\hss $\bu$\hss}}%
    \graphtemp=.5ex\advance\graphtemp by 0.680in
    \rlap{\kern 0.705in\lower\graphtemp\hbox to 0pt{\hss $\bu$\hss}}%
    \graphtemp=.5ex\advance\graphtemp by 0.680in
    \rlap{\kern 0.276in\lower\graphtemp\hbox to 0pt{\hss $\bu$\hss}}%
    \special{pn 11}%
    \special{pa 490 557}%
    \special{pa 490 309}%
    \special{fp}%
    \special{pa 490 557}%
    \special{pa 705 680}%
    \special{fp}%
    \special{pa 490 557}%
    \special{pa 276 680}%
    \special{fp}%
    \graphtemp=.5ex\advance\graphtemp by 1.051in
    \rlap{\kern 0.490in\lower\graphtemp\hbox to 0pt{\hss $\bu$\hss}}%
    \graphtemp=.5ex\advance\graphtemp by 0.309in
    \rlap{\kern 0.062in\lower\graphtemp\hbox to 0pt{\hss $\bu$\hss}}%
    \graphtemp=.5ex\advance\graphtemp by 0.309in
    \rlap{\kern 0.919in\lower\graphtemp\hbox to 0pt{\hss $\bu$\hss}}%
    \special{pa 490 1051}%
    \special{pa 62 309}%
    \special{fp}%
    \special{pa 62 309}%
    \special{pa 919 309}%
    \special{fp}%
    \special{pa 919 309}%
    \special{pa 490 1051}%
    \special{fp}%
    \graphtemp=.5ex\advance\graphtemp by 0.656in
    \rlap{\kern 0.490in\lower\graphtemp\hbox to 0pt{\hss $v$\hss}}%
    \graphtemp=.5ex\advance\graphtemp by 1.200in
    \rlap{\kern 0.490in\lower\graphtemp\hbox to 0pt{\hss $G'$\hss}}%
    \graphtemp=.5ex\advance\graphtemp by 0.557in
    \rlap{\kern 0.490in\lower\graphtemp\hbox to 0pt{\hss $\bu$\hss}}%
    \graphtemp=.5ex\advance\graphtemp by 0.309in
    \rlap{\kern 1.554in\lower\graphtemp\hbox to 0pt{\hss $\bu$\hss}}%
    \graphtemp=.5ex\advance\graphtemp by 0.680in
    \rlap{\kern 1.768in\lower\graphtemp\hbox to 0pt{\hss $\bu$\hss}}%
    \graphtemp=.5ex\advance\graphtemp by 0.680in
    \rlap{\kern 1.340in\lower\graphtemp\hbox to 0pt{\hss $\bu$\hss}}%
    \special{pa 1554 309}%
    \special{pa 1768 680}%
    \special{fp}%
    \special{pa 1768 680}%
    \special{pa 1340 680}%
    \special{fp}%
    \special{pa 1340 680}%
    \special{pa 1554 309}%
    \special{fp}%
    \graphtemp=.5ex\advance\graphtemp by 1.051in
    \rlap{\kern 1.554in\lower\graphtemp\hbox to 0pt{\hss $\bu$\hss}}%
    \graphtemp=.5ex\advance\graphtemp by 0.309in
    \rlap{\kern 1.126in\lower\graphtemp\hbox to 0pt{\hss $\bu$\hss}}%
    \graphtemp=.5ex\advance\graphtemp by 0.309in
    \rlap{\kern 1.983in\lower\graphtemp\hbox to 0pt{\hss $\bu$\hss}}%
    \special{pa 1554 1051}%
    \special{pa 1126 309}%
    \special{fp}%
    \special{pa 1126 309}%
    \special{pa 1983 309}%
    \special{fp}%
    \special{pa 1983 309}%
    \special{pa 1554 1051}%
    \special{fp}%
    \graphtemp=.5ex\advance\graphtemp by 0.557in
    \rlap{\kern 1.554in\lower\graphtemp\hbox to 0pt{\hss $T$\hss}}%
    \graphtemp=.5ex\advance\graphtemp by 1.200in
    \rlap{\kern 1.554in\lower\graphtemp\hbox to 0pt{\hss $G$\hss}}%
    \graphtemp=.5ex\advance\graphtemp by 0.309in
    \rlap{\kern 2.717in\lower\graphtemp\hbox to 0pt{\hss $\bu$\hss}}%
    \graphtemp=.5ex\advance\graphtemp by 0.680in
    \rlap{\kern 2.931in\lower\graphtemp\hbox to 0pt{\hss $\bu$\hss}}%
    \graphtemp=.5ex\advance\graphtemp by 0.680in
    \rlap{\kern 2.503in\lower\graphtemp\hbox to 0pt{\hss $\bu$\hss}}%
    \special{pa 2717 557}%
    \special{pa 2717 309}%
    \special{fp}%
    \special{pa 2717 557}%
    \special{pa 2931 680}%
    \special{fp}%
    \special{pa 2717 557}%
    \special{pa 2503 680}%
    \special{fp}%
    \graphtemp=.5ex\advance\graphtemp by 0.062in
    \rlap{\kern 2.717in\lower\graphtemp\hbox to 0pt{\hss $\bu$\hss}}%
    \graphtemp=.5ex\advance\graphtemp by 1.051in
    \rlap{\kern 2.717in\lower\graphtemp\hbox to 0pt{\hss $\bu$\hss}}%
    \graphtemp=.5ex\advance\graphtemp by 0.557in
    \rlap{\kern 2.717in\lower\graphtemp\hbox to 0pt{\hss $\bu$\hss}}%
    \special{pa 2717 309}%
    \special{pa 2717 62}%
    \special{fp}%
    \special{pa 2717 62}%
    \special{pa 2931 680}%
    \special{fp}%
    \special{pa 2931 680}%
    \special{pa 2717 1051}%
    \special{fp}%
    \special{pa 2717 1051}%
    \special{pa 2503 680}%
    \special{fp}%
    \special{pa 2503 680}%
    \special{pa 2717 62}%
    \special{fp}%
    \graphtemp=.5ex\advance\graphtemp by 0.656in
    \rlap{\kern 2.717in\lower\graphtemp\hbox to 0pt{\hss $v$\hss}}%
    \graphtemp=.5ex\advance\graphtemp by 1.200in
    \rlap{\kern 2.717in\lower\graphtemp\hbox to 0pt{\hss $G'$\hss}}%
    \graphtemp=.5ex\advance\graphtemp by 0.309in
    \rlap{\kern 3.459in\lower\graphtemp\hbox to 0pt{\hss $\bu$\hss}}%
    \graphtemp=.5ex\advance\graphtemp by 0.680in
    \rlap{\kern 3.673in\lower\graphtemp\hbox to 0pt{\hss $\bu$\hss}}%
    \graphtemp=.5ex\advance\graphtemp by 0.680in
    \rlap{\kern 3.245in\lower\graphtemp\hbox to 0pt{\hss $\bu$\hss}}%
    \special{pa 3459 309}%
    \special{pa 3673 680}%
    \special{fp}%
    \special{pa 3673 680}%
    \special{pa 3245 680}%
    \special{fp}%
    \special{pa 3245 680}%
    \special{pa 3459 309}%
    \special{fp}%
    \graphtemp=.5ex\advance\graphtemp by 0.062in
    \rlap{\kern 3.459in\lower\graphtemp\hbox to 0pt{\hss $\bu$\hss}}%
    \graphtemp=.5ex\advance\graphtemp by 1.051in
    \rlap{\kern 3.459in\lower\graphtemp\hbox to 0pt{\hss $\bu$\hss}}%
    \graphtemp=.5ex\advance\graphtemp by 0.557in
    \rlap{\kern 3.459in\lower\graphtemp\hbox to 0pt{\hss $T$\hss}}%
    \graphtemp=.5ex\advance\graphtemp by 1.200in
    \rlap{\kern 3.459in\lower\graphtemp\hbox to 0pt{\hss $G$\hss}}%
    \special{pa 3459 309}%
    \special{pa 3459 62}%
    \special{fp}%
    \special{pa 3459 62}%
    \special{pa 3673 680}%
    \special{fp}%
    \special{pa 3673 680}%
    \special{pa 3459 1051}%
    \special{fp}%
    \special{pa 3459 1051}%
    \special{pa 3245 680}%
    \special{fp}%
    \special{pa 3245 680}%
    \special{pa 3459 62}%
    \special{fp}%
    \graphtemp=.5ex\advance\graphtemp by 0.309in
    \rlap{\kern 4.449in\lower\graphtemp\hbox to 0pt{\hss $\bu$\hss}}%
    \graphtemp=.5ex\advance\graphtemp by 0.557in
    \rlap{\kern 4.201in\lower\graphtemp\hbox to 0pt{\hss $\bu$\hss}}%
    \graphtemp=.5ex\advance\graphtemp by 0.680in
    \rlap{\kern 4.449in\lower\graphtemp\hbox to 0pt{\hss $\bu$\hss}}%
    \graphtemp=.5ex\advance\graphtemp by 0.557in
    \rlap{\kern 4.696in\lower\graphtemp\hbox to 0pt{\hss $\bu$\hss}}%
    \graphtemp=.5ex\advance\graphtemp by 1.051in
    \rlap{\kern 4.449in\lower\graphtemp\hbox to 0pt{\hss $\bu$\hss}}%
    \special{pa 4449 309}%
    \special{pa 4201 557}%
    \special{fp}%
    \special{pa 4449 309}%
    \special{pa 4449 680}%
    \special{fp}%
    \special{pa 4449 309}%
    \special{pa 4696 557}%
    \special{fp}%
    \special{pa 4449 1051}%
    \special{pa 4201 557}%
    \special{fp}%
    \special{pa 4449 1051}%
    \special{pa 4449 680}%
    \special{fp}%
    \special{pa 4449 1051}%
    \special{pa 4696 557}%
    \special{fp}%
    \graphtemp=.5ex\advance\graphtemp by 0.210in
    \rlap{\kern 4.449in\lower\graphtemp\hbox to 0pt{\hss $v$\hss}}%
    \graphtemp=.5ex\advance\graphtemp by 1.200in
    \rlap{\kern 4.449in\lower\graphtemp\hbox to 0pt{\hss $G'$\hss}}%
    \graphtemp=.5ex\advance\graphtemp by 0.557in
    \rlap{\kern 4.943in\lower\graphtemp\hbox to 0pt{\hss $\bu$\hss}}%
    \graphtemp=.5ex\advance\graphtemp by 0.680in
    \rlap{\kern 5.191in\lower\graphtemp\hbox to 0pt{\hss $\bu$\hss}}%
    \graphtemp=.5ex\advance\graphtemp by 0.557in
    \rlap{\kern 5.438in\lower\graphtemp\hbox to 0pt{\hss $\bu$\hss}}%
    \graphtemp=.5ex\advance\graphtemp by 1.051in
    \rlap{\kern 5.191in\lower\graphtemp\hbox to 0pt{\hss $\bu$\hss}}%
    \special{pa 5191 1051}%
    \special{pa 4943 557}%
    \special{fp}%
    \special{pa 4943 557}%
    \special{pa 5191 680}%
    \special{fp}%
    \special{pa 5191 680}%
    \special{pa 5438 557}%
    \special{fp}%
    \special{pa 5438 557}%
    \special{pa 5191 1051}%
    \special{fp}%
    \special{pa 5191 1051}%
    \special{pa 5191 680}%
    \special{fp}%
    \graphtemp=.5ex\advance\graphtemp by 0.581in
    \rlap{\kern 5.191in\lower\graphtemp\hbox to 0pt{\hss $T$\hss}}%
    \graphtemp=.5ex\advance\graphtemp by 1.200in
    \rlap{\kern 5.191in\lower\graphtemp\hbox to 0pt{\hss $G$\hss}}%
    \special{pn 8}%
    \special{ar 5191 927 445 445 -2.159827 -0.981765}%
    \graphtemp=.5ex\advance\graphtemp by 1.324in
    \rlap{\kern 0.985in\lower\graphtemp\hbox to 0pt{\hss (a)\hss}}%
    \graphtemp=.5ex\advance\graphtemp by 1.324in
    \rlap{\kern 3.088in\lower\graphtemp\hbox to 0pt{\hss (b)\hss}}%
    \graphtemp=.5ex\advance\graphtemp by 1.324in
    \rlap{\kern 4.820in\lower\graphtemp\hbox to 0pt{\hss (c)\hss}}%
    \hbox{\vrule depth1.324in width0pt height 0pt}%
    \kern 5.500in
  }%
}%
}
\caption{Three cases for $G'$ and $G$ in
Lemma~\ref{plan-decomp}.\label{figplan}}
\end{figure}

Finally, the additional edges may be shared in pairs (Figure~\ref{figplan}c).
Now the component of $G'$ containing $v$ is $K_{2,3}$, and the component of $G$
containing $T$ is $K_4$.  Form $G''$ from $G$ by deleting this component; $G''$
has fewer facial triangles than $G$.  If $G''$ has at least three vertices,
then it decomposes into at most $2n-12$ edges and facial triangles, yielding a
decomposition of size at most $2n-8$ for $G$.  If $G''$ has at most two
vertices, then it is $K_1$, $2K_1$ or $K_2$; in each case, $G$ decomposes into
at most $2n-6$ edges and facial triangles.
\end{proof}

\begin{corollary}\label{mindeg3}
  If $G$ is an $n$-vertex plane graph with $\delta(G)\ge3$, then
  $\pol(G)\le2n-4$, with equality only if every facial cycle is
  a 4-cycle or $G = K_4$.  The same holds for $\theta_1(G)$.
\end{corollary}
\begin{proof}
Immediate from Lemma~\ref{plan-decomp} and the Decomposition Bound
(Lemma~\ref{clique-decomp}).
\end{proof}

In many cases, we will obtain bounds on $\ol(G)$ from bounds on $\pol(G)$.
The next several remarks and computations facilitate characterization of the
$n$-vertex planar graphs and the $n$-vertex graphs with largest pure overlap
number.

\begin{observation}\label{components}
  With the convention that $\pol(K_1)=1$, pure overlap number and intersection
  number are additive under disjoint union; in particular, $\pol(2K_2)=6$.
  Overlap number is not additive under disjoint union:  $\ol(K_2)=3$, but
  $\ol(2K_2)=5$.
\end{observation}

We observed in Corollary~\ref{edgebound} that $\pol(G)=\C{E(G)}$ when $G$ is
triangle-free and $\delta(G)=2$.  When $\delta(G)=1$, more labels may be needed.

\begin{proposition}\label{path}
  If $n \ge 2$, then $\pol(P_n)=n+1$.  If $n\ge3$, then $\pol(C_n)=n$.
\end{proposition}
\begin{proof} 
A representation using the sets $\{i,i+1\}$ for $1\le i\le n$ provides the
upper bound.  Since a \porep\ is an \irep\ and $P_n$ is triangle-free, each
label is used at most twice.  The endpoints of an edge have a common label used 
at no other vertex.  Also, each endpoint of $P_n$ has a label used nowhere else.
Hence $|f|$ must exceed the number of edges by $2$.

If $n\ge4$, then $C_n$ is triangle-free and $\delta(C_n)=2$, so
$\pol(C_n)=\C{E(G)}=n$.  For $C_3$, using the three $2$-sets in $\{1,2,3\}$
yields $\pol(C_n)=n$ again.
\end{proof}

\begin{proposition}\label{pol-star}
  If $m \ge 2$, then $\pol(K_{1,m}) = 2m$.  If $G-v=K_{1,m}$ and $d_G(v)\le 2$,
  then $\pol(G)\le 2m+1$.
\end{proposition}
\begin{proof}
Since a \porep\ is an \irep, the sets for leaves of $K_{1,m}$ are disjoint.
They must overlap the central set, so each has size at least $2$.  Equality
holds using two labels at each leaf and putting one label from each leaf in the
central set.

If $G-v=K_{1,m}$ and $G\ne K_{1,m+1}$, then $G$ is obtained from $C_3$, $C_4$,
or $P_4$ by appending leaves at one vertex of degree $2$.  By the Deletion
Bound (Lemma~\ref{poldeg2}), each leaf costs at most two new labels; combined
with Proposition~\ref{path}, this yields $\pol(G)\le 2m+1$.
\end{proof}

A graph $G$ is {\it $k$-degenerate} if every subgraph of $G$ has a vertex of
degree at most $k$.

\begin{theorem}\label{planar-pol}
  If $G$ is an $n$-vertex planar graph with $n \geq 3$, then $\pol(G)\le2n-2$,
  with equality if and only if $G\in\{K_{1,n-1},2K_2,K_2+K_1\}$.
  Furthermore, if $G$ is not 2-degenerate, then $\pol(G)\le 2n-4$, with
  equality when $\delta(G) \ge 2$ only if $G$ has $2n-4$ edges or is $K_4$.
\end{theorem}

\begin{proof}
Observation~\ref{components} and Propositions~\ref{path}--\ref{pol-star} take
care of the case $n=3$ and confirm equality for $K_{1,n-1}$.  This provides a
basis for induction.  Suppose that $n\ge4$.
\looseness -1
If $\delta(G)\ge3$, then Corollary~\ref{mindeg3} yields $\pol(G)\le 2n-4$, with
equality only when $G$ has $2n-4$ edges or is $K_4$.

Hence we may assume that $\delta(G)\le 2$.  Let $v$ be a vertex of minimum
degree.  By the Deletion Bound and the induction hypothesis, 
$\pol(G)\le\pol(G-v)+2\le 2n-2$.  Equality requires $\pol(G-v)=2n-4$ and
$d(v)\in\{1,2\}$ and $G-v\in\{K_{1,n-2},2K_2,K_2+K_1\}$.

If $G-v=K_2+K_1$, then $G\in\{2K_2,P_3+P_1,C_3+P_1,P_4\}$, and $\pol(G)$ is
$6,5,4,5$, respectively, using Observation~\ref{components} and
Proposition~\ref{path}.  If $G-v=2K_2$, then $G\in\{P_3+P_2,P_5,C_3+P_2\}$, and
$\pol(G)$ is $7,6,6$, respectively, using the same facts.
If $G-v=K_{1,n-2}$ and $G\ne K_{1,n-1}$, then Proposition~\ref{pol-star}
states that $\pol(G)\le 2n-3$.

For the final statement, suppose that $G$ is not 2-degenerate.  Now $G-v$ has a
subgraph with minimum degree at least 3.  Hence $G-v$ is not 2-degenerate, and
the bound improves to $\pol(G)\le \pol(G-v)+2\le 2n-4$, with equality only if
$\pol(G-v)=2n-6$.

When also $\delta(G)=2$, we have $G\ne K_4$ and need to show $|E(G)|=2n-4$.
If $\delta(G-v)\ge2$, then by the induction hypothesis, $G-v$ has $2n-6$ edges
or is $K_4$.  In the former case $|E(G)|=2n-4$, since
$d(v)=\delta(G)=2$.  In the latter case, $\pol(G)\le 5=2n-5$ by using the five
sets $123,14,24,34,15$.

Since $\delta(G)=2$ prohibits $\delta(G-v)=0$, the remaining case is
$\delta(G-v)=1$, and a leaf $u$ in $G-v$ must be a neighbor of $v$ in $G$.
Consider a \porep\ $f$ of $G-v$ using $2n-6$ elements; since $u$
is a leaf, $f(u)$ must have an element $a$ assigned to no other vertex.
Assign $a$ and a new element $b$ to $v$, and add $b$ to the set for the other
neighbor of $v$.  Since no assigned set contains another, this use of $b$
causes no trouble.  We have extended $f$ to a \porep\ of $G$ with only one new
label, so $\pol(G) = 2n-5$.   
\end{proof}

\begin{corollary}
If $G$ is a planar $n$-vertex graph with $n \geq 3$, then $\ol(G)\le 2n - 2$.
\end{corollary}

The upper bounds on pure overlap number simplify some cases in proving the
best upper bound on overlap number, which in general is smaller by 1.

\begin{theorem}\label{planar}
If $G$ is a planar $n$-vertex graph and $n\ge5$, then $\ol(G)\le 2n-5$, which
is sharp for $n=8$ and $n\ge 10$.
\end{theorem}
\begin{proof}
Example~\ref{planar-examp} establishes sharpness for $n=8$ and $n\ge10$ (and
$n=4$).  To prove the bound, we use induction on $n$, postponing the base case
to Proposition~\ref{small} below.  For $n>5$, we may assume $\delta(G)\ge 2$ by
Observation~\ref{freevertices}, Lemma~\ref{ol-leaves}, and the induction
hypothesis.

If $G$ is 2-degenerate, then $|E(G)|\le 2n-3$.  If $|E(G)|<2n-3$, then the Edge
Bound (Lemma~\ref{ol-edge-upper}) yields $\ol(G)\le2n-5$, so we may assume
equality.  By Euler's Formula, $G$ contains a triangle $T$.  If every vertex of
$T$ has a neighbor outside $T$, then each vertex of $G$ is incident to at least
two subgraphs in the decomposition of $G$ consisting of $T$ and $2n-6$
individual edges.  The Decomposition Bound (Lemma~\ref{clique-decomp}) now
yields $\pol(G)\le 2n-5$, and hence $\ol(G)\le2n-5$.

If $|E(G)|=2n-3$ and every triangle has a vertex $v$ with $d_G(v)=2$, then
$G-v$ has the same property, and by induction $G$ is the book $B_n$ and
$\ol(G)=3$ (Lemma~\ref{triangles}).

In the remaining case, $G$ is not 2-degenerate and $\delta(G) \ge 2$.
By Theorem~\ref{planar-pol}, either $\ol(G)\le \pol(G)\le 2n-5$,
or $|E(G)|=2n-4$ and the Edge Bound applies.    
\end{proof}

When the complement of a graph $G$ is edge-transitive, $G^+$ denotes the graph
obtained by adding any edge of the complement to $G$.

\begin{proposition}\label{small}
  If $G$ is an $n$-vertex graph, where $n\in\{4,5\}$, then $\ol(G)\le 2n-5$,
  except that $\ol(G)=4$ for $G\in \{P_4,K_4,K_{1,3}^+\}$.
\end{proposition}
\begin{proof}
If $G$ is a forest, then Theorem~\ref{tree} suffices.  Note also that
$\ol(K_3)=3$, and we may assume that $G$ has no isolated vertex or repeated
neighborhood, by Observation~\ref{freevertices}.

If $n=4$ and $G$ is not a forest and has no isolated vertex, then
$G\in\{C_4,C_4^+,K_{1,3}^+,K_4\}$.  Each graph has an edge, so
$\ol(G)\ge3$.  For $C_4$ and $C_4^+$, the repeated neighborhoods let three
elements suffice.  For $K_4$, we need an intersecting family of four
incomparable sets, which does not exist in $\{1,2,3\}$, but $\{123,41,42,43\}$
suffices.  For $K_{1,3}^+$, the triangle can only be represented in subsets of
$\{1,2,3\}$ using $\{12,23,13\}$, and no fourth subset intersects just one of
these.  Hence four elements are needed, and $\{123,124,13,23\}$ is an \orep.

For $n=5$, if $G$ has a vertex $v$ of degree 1 such that $\ol(G-v)\le3$, then
Lemma~\ref{ol-leaves} applies.  With no repeated neighborhood, a vertex of
degree 1 now restricts $G$ to be $K_4$ plus one pendant edge, $K_3$ plus
pendant edges at two distinct vertices, or $K_3$ plus a pendant path of length
two at one vertex.  These three graphs are represented by
$\{145,245,345,1234,45\}$, $\{12,23,34,45,1245\}$, and $\{12,23,34,45,1235\}$,
respectively.

We are left with $n=5$ and $\delta(G)\ge 2$.  By Lemma~\ref{ol-edge-upper}, we
may assume that $\C{E(G)}\ge7$.  The remaining 5-vertex graphs with at least
seven edges are listed below with \orep s (``$+$'' denotes disjoint union).

\medskip

\begin{tabular}{rlcrl}
$K_5$&$\{123,234,345,451,512\}$&&
    $K_{2,2,1}$&$\{12,34,14,23,13\}$\\
$\overline{P_2+3K_1}$&$\{123,234,345,14,25\}$&&
    $K_{3,1,1}$&$\{12,34,1234,513,524\}$\\
$\overline{P_3+2K_1}$&$\{123,345,14,25,1245\}$&&
    $\overline{P_4+K_1}$&$\{12,23,34,45,135\}$
  \end{tabular}
\end{proof}

\section{Extremal values for $n$-vertex graphs}

In this section we study the maximum values for $\pol(G)$ and $\ol(G)$ over
$n$-vertex graphs.  As usual, the problem is easier for $\pol(G)$, and solving
it simplifies the analysis for $\ol(G)$.  In addition to our earlier
computations, we need one more special family.

\begin{proposition}\label{clique}
  If $n\ge1$ and $\binom{2k-1}k\ge n$, then $\pol(K_n)\le 2k-1$.
\end{proposition}
\begin{proof}
The $k$-element subsets of a $(2k-1)$-set are pairwise intersecting and
incomparable.
\end{proof} 

\begin{lemma}\label{pol-triangle}
  Let $T$ be the vertex set of a triangle in an $n$-vertex graph $G$.  If $n'$
  is the number of vertices of degree $1$ in $G-T$, then
  $\pol(G)\le \pol(G')+n-n'$.
\end{lemma}
\begin{proof}
Let $f$ be an optimal \porep\ of $G-T$.  Add three new labels, two assigned to
each vertex, to represent the triangle.  Consider each vertex $x$ outside $T$.
If $d(x)=1$, then $f(x)$ has at least two labels, with one appearing on no
other vertex; add that label to the sets for the neighbors of $x$ in $T$.
If $d(x)\ne1$, then introduce a new label assigned to $x$ and its neighbors in
$T$.  The total number of labels used is at most $\pol(G-T)+3+n-3-n'$.
\end{proof}

\begin{theorem}\label{pol-upper-bound}
  Let $G$ be an $n$-vertex graph.  If $3\le n\le 5$, then $\pol(G)\le 2n-3$
  unless $G\in\{K_{1,n-1},2K_2,K_2+K_1\}$.  If $n \ge 6$ and $G\ne K_{1,5}$,
  then $\pol(G) \le \floor{n^2/4}$.  If $n\ge7$, then equality holds only
  when $G=K_{\floor{n/2},\ceil{n/2}}$.
\end{theorem}
\begin{proof} The first statement was proved in Theorem~\ref{planar-pol} except
for the nonplanar graph $K_5$, and $\pol(K_5)\le5$, by Proposition~\ref{clique}.
For $n\ge6$, we proceed inductively.

{\it Case 1: There exists $v\in V(G)$ with $d(v)\le 2$.}  By the Deletion Bound
(Lemma~\ref{poldeg2}), $\pol(G)\le\pol(G-v)+2$.  If $n=6$, then $\pol(G-v)\le7$
unless $G-v=K_{1,4}$; in either case, $\pol(G)\le9=6^2/4$
(Proposition~\ref{pol-star}).  If $n\ge7$, then
$\pol(G)\le\floor{(n-1)^2/4}+2<\floor{n^2/4}$ unless
$G-v=K_{1,5}$, in which case $\pol(G)\le 11<\floor{7^2/4}$
(Proposition~\ref{pol-star}).

{\it Case 2: $G$ is triangle-free and $\delta(G)\ge3$.}  By the Decomposition
Bound and Corollary~\ref{edgebound} and the well-known fact that
$K_{\floor{n/2},\ceil{n/2}}$ is the unique triangle-free $n$-vertex graph with
the most edges, $\pol(G)\le\floor{n^2/4}$, with equality only for
$K_{\floor{n/2},\ceil{n/2}}$.

{\it Case 3: $G$ has a triangle and $\delta(G) \ge 3$.}
Let $T$ be the vertex set of a triangle.  By Lemma~\ref{pol-triangle},
$\pol(G)\le\pol(G-T)+n-n'$, where $n'$ is the number of vertices of degree $1$
in $G-T$.  For $n\ge9$, we have $\pol(G)\le\floor{(n-3)^2/4}+n<\floor{n^2/4}$
unless $G-T=K_{1,5}$, in which case $n'=5$ and $\pol(G)\le 10+1=11$.

For $6\le n\le 8$, we have $\pol(G-T)+n\le 2n-9+n\le\floor{n^2/4}$ unless
$G-T\in\{K_{1,n-1},2K_2,K_2+K_1\}$.  In those cases, $\pol(G-T)=2n-8$ and
$n'\ge2$, so $\pol(G-T)<\floor{n^2/4}$.  In the remaining cases with $n=8$, we
have $\pol(G)\le 3n-9<\floor{n^2/4}$.

If $n=7$, then $\floor{n^2/4}=12$, so it remains only to prove $\pol(G)\le11$
when $\pol(G-T)=2n-9=5$ and $\delta(G-T)\ge2$.  From $\delta(G-T)=2$, we have
$G-T\in\{C_4,C_4^+,K_4\}$.  We have shown $\pol(C_4)=4$, and also
$\pol(C_4^+)=4$ by adding a set with one new label and one old label to the
\porep\ of $C_3$ using three labels.  Hence only $G-T=K_4$ remains.  Since this
must hold for every triangle in $G$, we have $G=K_7$, but $\pol(K_7)=5$.
\end{proof}

Our remaining task is to find the maximum of $\ol(G)$ over $n$-vertex graphs.
Rosgen~\cite{Ros} showed $\ol(G) \leq n^2/4$.  In Corollary~\ref{bicliq-match},
we constructed for $n\ge6$ an example having overlap number
$\floor{n^2/4-n/2-1}$.  We will improve the upper bound to show that it is
extremal for $n\ge14$.  We consider the main cases in separate lemmas:
bipartite graphs, triangle-free non-bipartite graphs, and graphs containing a
triangle.

\begin{lemma}\label{bipartite-distinct-nbrhd}
  Let $G$ be an $n$-vertex bipartite graph in which no two vertices
  have the same neighborhood. If $n \ge 7$ and $\delta(G) \ge 2$, then
  $\ol(G) \le \floor{n^2/4 - n/2 - 1}$.
\end{lemma}
\begin{proof}
By the Edge Bound (Lemma~\ref{ol-edge-upper}), $\ol(G)\le\C{E(G)}-1$, so we may
assume that $\C{E(G)}>\floor{n^2/4-n/2}$.  Let $X$ and $Y$ be the partite sets,
with $k=\C X\le\C Y$.  To avoid duplicate neighborhoods, at most one vertex of
$Y$ has degree $k$, and vertices of degree $k-1$ have distinct nonneighbors
in $X$.  Summing the vertex degrees in $Y$ yields
$$
\C{E(G)}\le (k-1)(n-k)+1\le \floor{(n-1)^2/4}+1.
$$
Equality holds only when $k=\ceil{n/2}$, but we have $k\le \floor{n/2}$,
so we may assume $k=n/2$ and $n\ge8$.  Furthermore, $G$ arises from
$K_{n/2,n/2}$ by deleting a matching of size $n/2-1$.

Let $y$ be the vertex in $Y$ having degree $n/2$, and let $G'=G-y$.  Since
$n\ge8$, we have $\delta(G')\ge2$.  By the Decomposition Bound
(Lemma~\ref{clique-decomp}), $\pol(G')\le \C{E(G')}= n^2/4-n+1$.  Choose
$y'\in Y-\{y\}$, and let $x'$ be the nonneighbor of $y'$ in $X$.  Let $f'$ be a
\porep\ of $G'$ using one label for each edge.  Define $f$ as follows: Put
$f(y)=f'(y')\cup\{a\}$, where $a$ is a new label, and let
$f(x')=f'(x')\cup\{a\}$.  For $v\notin\{y,x'\}$, let $f(v)=f'(v)$.

Since the only vertex of $G'$ receiving $a$ is $x'$, overlaps and disjointness
are preserved within the sets assigned to $V(G')$ Hence it suffices to check
pairs involving $y$.  We have $f(y) \supset f(y')$ and $f(y) \disj f(z)$ for
$z \in Y-\{y'\}$.  For $x \in X$, we have $f(y) \olap f(x)$.  Thus, $f$ is an
\orep\ of $G$ with $n^2/4-n+2$ labels.  Since $n\ge 8$, the desired bound holds.
\end{proof}

\begin{lemma}\label{ol-bipartite}
  If $G$ is an $n$-vertex bipartite graph, then
  $\ol(G)\le\max\{2n,\floor{n^2/4 - n/2 - 1}\}$.
\end{lemma}

\begin{proof} 
The claim is $\ol(G)\leq2n$ for $n\le10$ and $\ol(G)\le\floor{n^2/4-n/2-1}$ for
$n>10$.  Since $\ol(G)\le\pol(G)$, Theorem~\ref{pol-upper-bound} implies the
claim when $n\le 8$, using $2n-2\le 2n$ always and $\floor{n^2/4}\le 2n$ for
$n\le 8$.

We proceed inductively.  The desired bound exceeds the desired bound for
$(n-1)$-vertex graphs by at least $2$.  Consider $v,w\in V(G)$.  If $d(v)=0$,
then $\ol(G)\le \ol(G-v)+1$.  If $d(v)=1$, then $\ol(G)\le \ol(G-v)+2$ (by
Lemma~\ref{ol-leaves}).  If $N(v)=N(w)$, then Observation~\ref{freevertices}
yields $\ol(G)=\ol(G-v)$.  If $\delta(G) \ge 2$ and no two vertices have the
same neighborhood, then Lemma~\ref{bipartite-distinct-nbrhd} yields
$\ol(G)\le\floor{n^2/4-n/2-1}$.  Thus the bound holds in all cases.
\end{proof}

\begin{lemma}\label{ol-nonbipartite}
  If $G$ is an $n$-vertex triangle-free graph that is not bipartite, then
  $\ol(G) \le \max\{2n+7, \floor{n^2/4 - n/2 - 1}\}$.
\end{lemma}

\begin{proof}
The claim is $\ol(G)\leq2n+7$ for $n\le12$ and $\ol(G)\le\floor{n^2/4-n/2-1}$
for $n>12$.  As above, Theorem~\ref{pol-upper-bound} implies the claim when
$n\le 10$.  As in Lemma~\ref{ol-bipartite}, we proceed inductively; the desired
bound increases by at least $2$ per step, and we may assume that
$\delta(G)\ge2$ and that no two vertices have the same neighborhood.  By
the Edge Bound (Lemma~\ref{ol-edge-upper}), $\ol(G)\le\C{E(G)}-1$, so it
suffices to show $\C{E(G)}\le n^2/4 - n/2$ for a triangle-free graph $G$ with
no repeated neighborhood.

Let $C$ be a shortest odd cycle in $G$, with length $2k+1$, and let $G'=G-V(C)$.
Since $C$ has no chords, $V(C)$ induces $2k+1$ edges.  Since $G$ has no
triangle, each vertex not on $C$ has at most $k$ neighbors on $C$.  Since
$G'$ is triangle-free, $\C{E(G')}\le(n-2k-1)^2/4$.  Summing the bounds $2k+1$,
$k(n-2k-1)$, and $(n-2k-1)^2/4$ yields $\C{E(G)}\le n^2/4-n/2-(k^2-2k-5/4)$.
If $k\ge 3$, then $k^2-2k>5/4$, so we may assume $k=2$.

If $G'$ is not bipartite, then let $C'$ be a shortest odd cycle in $G'$, with
length $2l+1$.  With $\C{V(G')-V(C')}=n-2l-6$, we have
$\C{E(G)}\le 5 + 2(n-5) + (2l+1) + l(n-2l-6) + (n-2l-6)^2/4$.  The bound
simplifies to $n^2/4 - n+5-l(l-2)$.  Since $l\ge2$, and $n-5\ge n/2$ when
$n\ge10$, it is small enough.

Finally, suppose that $G'$ is bipartite.  Since $k=2$, we have
$\C{E(G)}\le n^2/4-n/2+5/4$.  Call a vertex of $G'$ {\it full} if it has two
(nonadjacent) neighbors in $C$ and is adjacent to all vertices in the other
partite set of $G'$.  Each pair of nonadjacent vertices in $C$ is adjacent to
at most one full vertex, since otherwise $G$ has a triangle or a repeated
neighborhood.  Thus, at most five vertices of $G'$ are full, so at least
$(n-5)-5$ vertices are not.  This yields
$\C{E(G)}\le \floor{n^2/4-n/2+5/4-(n-10)/2}$, which suffices when $n\ge11$.
\end{proof}

\begin{lemma}\label{ol-triangle}
  If $G$ is an $n$-vertex graph with a triangle $T$, and $n\ge14$, then
  $\ol(G) \le \floor{n^2/4 - n/2 - 1}$.
\end{lemma}
\begin{proof} 
View $T$ as a triple of pairwise adjacent vertices, and let $G'=G-T$.
In several cases we show $\pol(G)\le\floor{n^2/4-n/2-1}$, which suffices.

{\it Case 1: $G'$ has a triangle $T'$.} By Theorem~\ref{pol-upper-bound},
$\pol(G'-T')\le\floor{(n-6)^2/4}$ when $n\ge13$.  By Lemma~\ref{pol-triangle}
(twice), $\pol(G)\le \pol(G'-T')+2n-3$.  If $n\ge 14$, then
${(n-6)^2/4}+2n-3 \leq {n^2/4-n/2-1}$.

{\it Case 2: $\delta(G')\le1$.} 
By Theorem~\ref{pol-upper-bound}, $\pol(G'-v)\le\floor{(n-4)^2/4}$ when
$n\ge11$.  By the Deletion Bound and then Lemma~\ref{pol-triangle},
$\pol(G)\le\pol(G'-v)+n+1$.  Hence $\pol(G)\le n^2/4-n+5$.  Since
$n-5\ge n/2+1$ when $n\ge12$, we have $\pol(G)\le n^2/4-n/2-1$.

{\it Case 3: $\delta(G')\ge2$ and $G'$ is triangle-free but not bipartite.} 
The argument in the second paragraph of Lemma~\ref{ol-nonbipartite} yields
$\C{E(G')}\le (n-3)^2/4 - (n-3)/2+5/4 = n^2/4 - 2n + 5$.  By the Decomposition
Bound (Lemma~\ref{clique-decomp}, $\pol(G')\le\C{E(G')}$.  By
Lemma~\ref{pol-triangle}, $\pol(G)\le \pol(G')+n\le  n^2/4 - n + 5$.  As in
Case 2, the claim holds when $n\ge 12$.

{\it Case 4: $\delta(G')\ge2$ and $G'$ is bipartite.} 
Suppose first that $T\esub N_G(v)$ for some $v\in V(G')$.  Since $K_4$ has a
\porep\ using $\{123,41,42,43\}$, the method of Lemma~\ref{pol-triangle} yields
$\pol(G)\le \pol(G'-v)+n\le (n-4)^2/4+n=n^2/4-n+4$.  This suffices when
$n\ge11$.

Thus, we may assume that each vertex of $G'$ has at most two neighbors in $T$.
Using Lemma~\ref{pol-triangle}, our present bound on $\pol(G)$ is 
$\C{E(G')}+n$, and $\C{E(G')}\le (n-3)^2/4$, so
$\pol(G)\le n^2/4-n/2+9/4$, and we only need to reduce this by $13/4$.
 
Call a vertex of $G'$ {\it full} if it has at least one neighbor in $T$ and is
adjacent to all vertices in the other partite set of $G'$.  Each vertex that is
not full reduces the added number of labels in the construction of
Lemma~\ref{pol-triangle} by 1 or reduces the degree-sum in $G'$ by $1$.  Hence
when $n\ge13$ it suffices to show that there are at most six full vertices.
We prove this for $n\ge14$.

Since the neighborhood in $T$ of a full vertex has size 1 or 2, there are only
six possible such neighborhoods.  Nonadjacent full vertices with the same
neighborhood in $T$ would have the same neighborhood in $G$, in which case
Observation~\ref{freevertices} completes the proof (using the inductive bound
on $\ol(G)$, not on $\pol(G)$).  Hence having seven full vertices requires full
vertices $u$ and $v$ that are adjacent in $G'$ and have the same neighborhood
$S$ in $T$.  We argue that this leads to two disjoint triangles in $G$, which
allows Case 1 to complete the proof.

Since nonadjacent full vertices cannot have the same neighborhood in $T$, no
other full vertex has neighborhood $S$ in $T$.  If two other adjacent full
vertices $x$ and $y$ have neighborhood $S'$ in $T$, then $S\ne S'$, so there
exist distinct vertices $s\in S$ and $s'\in S'$, and $\{s,u,v\}$ and
$\{s',x,y\}$ are disjoint triangles.

Otherwise, for each $S'\esub T$ with $|S'|\le 2$, some full vertex has
neighborhood $S'$ in $T$.  For $s\in S$, the triangle $\{s,u,v\}$ is disjoint
from the triangle formed by $T-\{s\}$ and the full vertex having neighborhood
$T-\{s\}$ in $T$.
\end{proof}

\begin{theorem}
If $n(G) \geq 14$, then $\ol(G)\leq n^2/4-n/2-1$.
\end{theorem}

\begin{proof} 
The claim follows immediately from Lemmas \ref{ol-bipartite},
\ref{ol-nonbipartite}, and \ref{ol-triangle}.
\end{proof}

We believe that the bound in fact holds for $n\ge8$ and is sharp only for the
construction in Corollary~\ref{bicliq-match} (deleting a perfect matching or a
near-perfect matching plus one edge from $K_{\floor{n/2},\ceil{n/2}}$).
Proving this seems likely to require substantial case analysis.


\end{document}